\RequirePackage{etoolbox}
\csdef{input@path}{{style/}{graphics/}}
\documentclass[MSNbibl,number,citesort,dvips]{arxbj}

\aid{0}
\volume{21}
\issue{1}
\pubyear{2015}
\firstpage{374}
\lastpage{400}
\doi{10.3150/13-BEJ571} 

\makeatletter

\newtheorem{theo}{Theorem}[section]
\newtheorem{lem}[theo]{Lemma}
\newremark{rem}[theo]{Remark}
\newtheorem{prop}[theo]{Proposition}
\newtheorem{cor}[theo]{Corollary}
\newproclaim{defn}{Definition}[section]

\renewcommand{\H}{{\mathcal H}}
\newcommand{\K}{{\mathcal K}}
\renewcommand{\S}{{\mathcal S}}
\renewcommand{\SS}{{\mathbb S}}
\renewcommand{\P}{{\mathbb P}}
\newcommand{\F}{{\mathcal F}}
\newcommand{\G}{{\mathcal G}}
\newcommand{\D}{{\mathcal D}}
\newcommand{\M}{{\mathcal M}}
\newcommand{\B}{\mathbf{B}}
\newcommand{\E}{{\mathbb E}}
\newcommand{\as}{\mbox{$\P$-a.s.}}
\newcommand{\R}{{\mathbb R}}
\newcommand{\X}{{\mathcal X}}
\newcommand{\C}{{\mathcal C}}
\newcommand{\I}{\mathbf{1}}
\newcommand{\Z}{{\mathbb{Z}}}

\newcommand{\eqref}[1]{(\ref{#1})}
\newcommand{\Mon}{\operatorname{Mon}}
\newcommand{\Vect}{\operatorname{Vect}}
\newcommand{\Ker}{\operatorname{Ker}}
\makeatother

\begin{document}
\begin{frontmatter}

\title{A compact LIL for martingales in $2$-smooth Banach spaces with
applications}
\runtitle{A compact LIL for Banach-valued martingales}

\begin{aug}
\author{\inits{C.}\fnms{Christophe} \snm{Cuny}\corref{}\ead[label=e1]{christophe.cuny@ecp.fr}}
\address{Laboratoire MAS, Ecole Centrale de Paris, Grande Voie des
Vignes, 92295 Chatenay-Malabry cedex, France. \printead{e1}}
\end{aug}

\received{\smonth{10} \syear{2012}}
\revised{\smonth{6} \syear{2013}}

%
\begin{abstract}
We prove the compact law of the iterated logarithm
for stationary and ergodic differences of (reverse or not) martingales
taking values in a separable $2$-smooth Banach space (for instance a
Hilbert space). Then, in the martingale case, the almost sure
invariance principle is
derived from a result of Berger. From those results, we deduce
the almost sure invariance principle for stationary processes under
the Hannan condition and the compact law of the iterated logarithm
for stationary processes arising from non-invertible dynamical systems.
Those results for stationary processes are new, even in the real valued
case. We also obtain the Marcinkiewicz--Zygmund strong law of large numbers
for stationary processes with values in some smooth Banach spaces.
Applications to several situations are given.
\end{abstract}

%

\end{frontmatter}

\section{Introduction}
Let $(\X,|\cdot|_\X)$ be a separable Banach
space and $\X^*$ be its
topological dual. Let $(\Omega,\F,\P)$ be a probability space and let
$(X_n)_{n\ge0}$ be a strictly stationary sequence of $\X$-valued random
variables. We are interested in the $\P$-a.s. behaviour of
$(S_n/\sqrt{2n L(L(n))})_{n\ge1}$, where $S_n:=X_0+\cdots+X_{n-1}$ and
$L:=\max(1,\log)$.

\begin{defn}We say that $(X_n)_{n\ge0}$ satisfies the
bounded law of the iterated logarithm (bounded LIL or BLIL) if
$(S_n/\sqrt{2n L(L(n))})_{n\ge1}$ is $\P$-a.s. bounded.
\end{defn}

\begin{defn}We say
that $(X_n)_{n\ge0}$ satisfies the
compact law of the iterated logarithm (compact LIL or CLIL) if
$(S_n/\sqrt{2n L(L(n))})_{n\ge1}$ is $\P$-a.s. relatively compact.
\end{defn}

When $(X_n)_{n\ge0}$ is a sequence of independent random variables,
the bounded and compact LILs are well understood, thanks to a characterization
due to Ledoux and Talagrand \cite{LT}. When the compact LIL holds, the
cluster set of $S_n/\sqrt{2n L (L(n))})_{n\ge1}$ may be identified
thanks to a result of Kuelbs \cite{Kuelbs}. When $X_0$ is \emph{pregaussian}
(see next section), we have an almost sure invariance principle as well.

For Banach spaces of type 2 (see next section for the definition),
the result of Ledoux--Talagrand takes the following particularly simple form.

\begin{theo}[(Ledoux and Talagrand {\cite{LTbook}, Corollary~8.8})]\label{LTtheo}
Let $(X_n)_{n\ge0}$ be a sequence of i.i.d. random variables with values
in a Banach space of type 2. Then, $(X_n)_{n\ge0}$ satisfies the
bounded LIL (resp. the
compact LIL) if and only if $\E((x^*(X_0))^2)<\infty$
for every $x^*\in\X^*$ (resp. $ ((x^*(X_0))^2)_{x^*\in\X
^*,|x^*|_{\X^*}\le1}$
is uniformly integrable), $\E(|X_0|_\X^2/L(L(|X_0|_\X)))<\infty$ and
\mbox{$\E(X_0)=0$}.
\end{theo}

In particular, a sequence of i.i.d. variables $(X_n)_{n\ge0}$ with
values in a Banach space of type 2
satisfies the compact LIL (hence, the bounded LIL) as soon as:
%
\begin{equation}
\label{mincondLIL} \E\bigl(|X_0|_\X^2\bigr)<
\infty\quad \mbox{and}\quad  \E(X_0)=0 .
\end{equation}

Now (see Remark~\ref{rempisier}), by a result of Pisier \cite
{Pisier1}, if $\X$ is a Banach space for which any sequence of
$\X$-valued i.i.d. variables, such that \eqref{mincondLIL} holds,
satisfies the bounded LIL, then, $\X$ must be of type
$p$ for any $1<p<2$.

We are interested here in the case where $(X_n)_{n\ge0}$ is a general
stationary sequence, including the case of martingale differences (and of
reverse martingale differences). The analogue of the notion of Banach space
of type 2
in the case of martingale differences is the
notion of $2$-smooth Banach space (see the next section for the definition).
One could wonder whether Theorem~\ref{LTtheo} is true in this context, or,
at least, whether \eqref{mincondLIL} is sufficient for the bounded LIL
or the
compact LIL, when $(X_n)_{n\ge0}$ is
a stationary sequence of martingale differences.

As far as we know, the latter question remained unsolved. Let us
however mention some results in that direction. Morrow and Philipp
\cite{MoP} (see also \cite{Philipp} for an improved version) obtained
an almost sure
invariance principle (see the next section for the definition), hence
a compact LIL (with an ad hoc normalization),
for sequences of non-necessarily stationary martingale differences taking
values in a Hilbert space. Dehling, Denker and Philipp \cite{DDP}
proved a
bounded LIL in the same context. When applied to
stationary sequences of martingale differences, the above results
require higher moments than 2.

In this paper, we prove that condition \eqref{mincondLIL}
is sufficient for the compact LIL when $(X_n)_{n\ge0}$ is a stationary sequence
of martingale
differences with values in a $2$-smooth Banach space. When the sequence is
ergodic, the cluster set of $(S_n/\sqrt{2n L(L(n))})_{n\ge1}$ is
identified as well as $\limsup_n|S_n|_\X/\sqrt{nL(L(n))}$. Then, using
a result of Berger \cite{Berger}, we obtain an almost sure invariance
principle for $(S_n)_{n\ge1}$. Those results (except for the
invariance principle) extend to \emph{reverse} martingale differences.
We do not know whether our results could be extended beyond
the scope of $2$-smooth Banach spaces. However, the above mentioned
result of
Pisier shows some limitation.

To prove those results, we first obtain integrability properties of
the ``natural" maximal function arising in that context, hence
generalizing a result of Pisier \cite{Pisier1} for i.i.d. variables.
This step is
crucial not only to prove the results for martingales
(and reverse martingales),
but also in order to extend the results to general stationary
processes under projective conditions, such as the Hannan condition, see
Theorem~\ref{theoHannan} or the Maxwell--Woodroofe condition, see
Cuny \cite{CunyMW}.
We note that the almost sure
invariance principle for Hilbert-valued stationary processes under
mixing conditions have been obtained by Merlev\`ede \cite{Merlevede}
and Dedecker and Merlev\`ede \cite{DMASIP}. Their results have
different range of
applications.

We also investigate the Marcinkiewicz--Zygmund strong law of large
numbers for stationary processes taking values in a smooth Banach
space. The maximal function arising in that
other context has been studied by
Woyczy\'nski \cite{Woyczinski}, for stationary martingale differences. We
investigate the case of stationary processes under projective conditions.
The main argument used is the
same as the one for the law of the iterated logarithm. The
Marcinkiewicz--Zygmund strong laws in smooth Banach spaces have been
also studied by Dedecker and Merlev\`ede \cite{DM} for stationary
processes satisfying
mixing conditions.

In the next section, we set our notations and state our results for
martingales and then, for stationary processes, including non-adapted
processes, functionals of Markov chains or
iterates of non-invertible dynamical systems. In Section~\ref{examples}, we give several examples to
which our conditions apply. In Section~\ref{proofs1}, we prove
our martingale results and in Section~\ref{proofs2} we prove our
results for stationary processes.
Finally, we postpone some technical proofs or results to the \hyperref[appA]{Appendix}.

\section{Main results}\label{results}

Let $(\Omega,\F,\P)$ be a probability space. We will consider
Banach-valued random variables. We refer to the book by Ledoux and
Talagrand \cite{LTbook} for the basic facts on the topic (definition,
conditional expectation\ldots).

Let $(\X,|\cdot|_\X)$ be a separable \emph{real} Banach space. We
endow $\X$ with its Borel $\sigma$-algebra. Denote by $L^0(\X)$ the
space (of classes modulo $\P$) of measurable random variables on
$\Omega$ taking values in $\X$. We define, for every $p\ge1$, the
usual Bochner spaces $L^p$ and their weak versions, as follows
\begin{eqnarray*}
L^p(\Omega,\F,\P,\X)&=&\bigl\{ Z\in L^0(\X) \dvt  \E
\bigl(|Z|_\X^p \bigr)<\infty\bigr\} ;
\\
L^{p,\infty}(\Omega,\F,\P,\X)&=&\Bigl\{ Z\in L^0(\X) \dvt  \sup
_{t >0} t\bigl(\P \bigl(|Z|_\X>t\bigr)\bigr)^{1/p} <
\infty\Bigr\} .
\end{eqnarray*}

For every $Z\in L^p(\Omega,\F,\P,\X)$, write $\|Z\|_{p,\X}:=(\E
(|Z|_\X^p ))^{1/p}$ and for
every $Z\in L^{p,\infty}(\Omega,\allowbreak \F,\P,\X)$, write $\|Z\|
_{p,\infty,\X}:= \sup_{t >0} t(\P(|Z|_\X>t))^{1/p}$.

For the sake of clarity, when they are understood, some of the references
to $\Omega$, $\F$ or $\P$ may be omitted. Also, in the case
when $\X=\R$, we shall simply write $\|\cdot\|_{p}$ or $\|\cdot\|_{p,
\infty}$.
Recall that for every $p>1$ there
exists a norm on $L^{p,\infty}(\P,\X)$ (see, for instance, \cite{LTbook},
Chapter ``Notation''), equivalent to the quasi-norm
$\|\cdot\|_{p,\infty,\X}$, that makes $L^{p,\infty}(\P,\X)$ a
Banach space.


The Banach spaces we will consider are the so-called \emph{smooth}
Banach spaces.

\begin{defn}
We say that $\X$ is $r$-\emph{smooth}, for some $1<r\le2$, if there
exists $L\ge1$,
such that
\[
|x+y|_\X^r+|x-y|_\X^r \le2
\bigl(|x|_\X^r+L^r|y|_\X^r
\bigr) \qquad \forall x,y\in \X .
\]
\end{defn}

\begin{defn}
We say that $(d_n)_{1\le n\le N}\subset L^1(\Omega,\F,\P,\X)$ is a
sequence of martingale differences, if there exists non-decreasing
$\sigma$-algebras $(\G_n)_{0\le n\le N}$ such that for every $1\le
n\le N$, $d_n$ is
$\G_n$-measurable and $\E(d_n|\G_{n-1})=0\ \as$ If $(\G_n)_{1\le
n\le N+1}$
is non-increasing and $\E(d_n|\F_{n+1})=0\ \as$, we speak about
differences of
\emph{reverse} martingales.
\end{defn}

It is known, see, for instance, Proposition~1 of Assouad \cite{Assouad}
(and its corollary), that when $\X$ is $r$-smooth, there exists $D\ge1$,
such that for every martingale differences $(d_n)_{1\le n\le N}$, we have
%
\begin{equation}
\label{smooth} \E\bigl(|d_1+\cdots+d_N|_\X^r
\bigr)\le D^r \sum_{n=1}^N\E
\bigl(|d_n|_\X^r\bigr) .
\end{equation}
When needed, we will say that $\X$ is $(r,D)$-smooth, where $D$ is a
constant such that condition \eqref{smooth} is satisfied (notice that
this definition is
compatible with the definition page 1680 of \cite{Pinelis}, see
Proposition~2.5
there). Clearly, $D$ must be greater than 1.


Any $L^p$ space, $p>1$ (of $\mathbb{R}$-valued functions), associated
with a $\sigma$-finite measure is $r$-smooth for $r=\min(2,p)$ (one
may take
$D^2=p-1$ if $p\ge2$, see \cite{Pinelis}, Proposition~2.1, and
$D^2=2$ if $1 \le p< 2$ by \cite{Assouad}). Any Hilbert space is
$(2,1)$-smooth. 

\begin{defn}
We say that $\X$ is a Banach space of type $r$, $1<r\le2$, if \eqref{smooth}
holds for every finite set $(d_n)_{1\le n\le N}$ of independent variables.
Hence, $2$-smooth Banach spaces are
particular examples of spaces of type $2$.
\end{defn}

Our goal is to study the law of the iterated logarithm and
the Marcinkiewicz--Zygmund strong law of large numbers
for the partial sums of an $\X$-valued stationary process. We will
start by studying the maximal functions associated
with these limit theorems. Let us specify some notations.

%


Let $\theta$ be a measurable measure preserving
transformation on $\Omega$. To any $X\in L^0(\Omega,\X)$, we
associate a stationary process $(X \circ\theta^n)_{n\ge0}$ (when
$\theta$ is invertible, we extend that definition to $n\in\mathbb
{Z}$). Then, for every $n\ge1$, write $S_n(X)=\sum_{i=0}^{n-1}
X\circ\theta^i$.

We shall assume that there exists a suitable filtration on $\Omega$.
In order to cover more situations, we shall consider filtrations that are
either non-decreasing or non-increasing. In spirit, the first case
arise when $\theta$ is invertible and the second one when $\theta$ is
non-invertible.

In particular, we assume that we are in one of the following situations.

If $\F_0\subset\F$ is a $\sigma$-algebra
such that $\F_0\subset\theta^{-1}(\F_0)$, we define a \emph
{non-decreasing} filtration
$(\F_n)_{n\ge0}$ by $\F_n:=\theta^{-n}(\F_0)$. Define then $\E_n=
\E(\cdot| \F_n)$.

If $\F^0$ is such that $\theta^{-1}(\F^0) \subset\F^0$ (for instance,
take $\F^0=\F$), we define a \emph{non-increasing} filtration $(\F
^n)_{n\ge0}$, by $\F^n:=\theta^{-n}(\F^0)$.
Define then $\E^n=\E(\cdot| \F^n)$.

Let $1 \le p \le2$. Let $X\in L^p(\Omega,\F,\P,\X)$. We consider the
following maximal functions
%
\begin{eqnarray}
\label{maxfunp} \M_{p}(X)&:=& \sup_{n\ge1}
\frac{| \sum_{k=0}^{n-1}
X\circ\theta^k|_\X}{n^{1/p}} , \qquad \mbox{if $1\le p<2$} ,
\\
\label{maxfun2} \M_{2}(X)&:=& \sup_{n\ge1}
\frac{| \sum_{k=0}^{n-1}
X\circ\theta^k|_\X}{\sqrt{n L(L(n))}} ,
\end{eqnarray}
where $L:=\max(\log,1)$.



The maximal operator $\M_1$ is related to Birkhoff's ergodic theorem, which
asserts that for every $X\in L^1(\Omega,\X)$, $((\sum_{k=0}^{n-1}X\circ
\theta^k)/n)_{n\ge1}$
converges $\P$-a.s. (see Theorem~2.1, page 167 of \cite{Krengel}
for the $\X$-valued case). For every $X\in L^1(\Omega,\X)$, by
Hopf's dominated ergodic theorem for real-valued stationary processes
(see \cite{Krengel}, Corollary 2.2, page 8),
applied to $(|X|_\X\circ\theta^n)_{n\ge0}$, we have
%
\begin{equation}
\label{Hopf} \bigl\|\M_1(X)\bigr\|_{1,\infty} \le\|X\|_{1,\X} .
\end{equation}
Now, once we know that \eqref{Hopf} holds, by the Banach principle (see
\cite{Krengel}, Theorem~7.2, page 64, or Proposition~\ref{banprin}), in
order to
prove Birkhoff's ergodic theorem, it suffices to
prove it on a set of $X$'s dense in $L^1$ (e.g., the $\theta$
invariant elements and the coboundaries). We want to use that strategy
to study the Marcinkiewicz--Zygmund strong law of large numbers and
versions of the law of the iterated logarithm. Of course, one cannot
expect to have a version of \eqref{Hopf} for ${\mathcal M}_p$, when
$1<p\le 2$
 without any further assumption on $(X\circ\theta
^n)_{n\ge0}$.




\subsection{Results for stationary (reverse) martingale differences}

In this subsection, we consider stationary sequences of (reverse) martingale
differences.

Let $d\in L^p(\Omega,\F_1,\X)$ be such that $\E_0(d)=0\ \as$ Then,
by our assumptions on $\F_0$, $(d\circ\theta^n)_{n\ge0}$ is a
stationary sequence of martingale differences.

Let $d\in L^p(\Omega,\F^0,\X)$ be such that $\E^1(d)=0\ \as$
Then, by
our assumption on $\F^0$, $(d\circ\theta^n)_{n\ge0}$ is a
stationary sequence of reverse martingale differences, that is, for
every $n\ge0$,
$d\circ\theta^n$ is $\F^n$-measurable and $\E(d\circ\theta^n|\F
^{n+1})=0\ \as$

There is no loss of generality in assuming that our stationary sequences
of (reverse) martingale differences are given that way.


Indeed, it is well known (see, e.g., Doob \cite{Doob}, page 456) that,
given a stationary sequence $(\tilde d_n)_{n\ge1}$ on a probability space
$(\tilde\Omega,\tilde\F,\tilde\P)$, there exist another
probability space
$(\Omega,\F,\P)$, an invertible bi-measurable
measure-preserving transformation $\theta$ on $ \Omega$ and a
random variable $d$ on $ \Omega$ such that the sequences
$(\tilde d_n)_{n\ge1}$ and $( d \circ\theta^n)_{n\ge1}$ have the same
law.

Moreover, it follows from the construction, that if $(\tilde d_n)_{n\ge
1}$ are
martingale
differences (respectively, reverse martingale differences), $( d \circ
\theta^n)_{n\ge1}$ are martingale
differences (respectively, reverse martingale differences) either.

Hence, since all the results we are concerned with in that paper only
rely on
the distribution of the processes under consideration, we shall assume
(without loss of generality) that
our stationary sequences of martingale differences are given thanks
to a measure-preserving transformation.

We start with a result of Woyczy\'nski about the Marcinkiewicz--Zygmund
strong law of large numbers.

\begin{prop}[(Woyczy\'nski \cite{Woyczinski})]\label{mzmax}
Let $1< p<r\le2$ and $D\ge1$. Let $\X$ be a separable $(r,
D)$-smooth Banach
space. There exists $C_{p,r}>0$ such that for every $d\in L^p(\Omega
,\F_1,\X)$ (resp. $d\in L^p(\Omega,\F^0,\X)$),
with $\E_0(d)=0$ (resp. $\E^1(d)=0$), we have
%
\begin{equation}
\label{weakMZ} \bigl\|\M_p(d)\bigr\|_{p,\infty}\le C_{p,r}
D^{r/p}\|d\|_{p,\X} .
\end{equation}
Moreover,
%
\begin{equation}
\label{MZ}\bigl|S_n(d)\bigr|_\X/n^{1/p}\to0 \qquad \as
\end{equation}
\end{prop}

\begin{rem} 
Actually, Woyczy\'nski proved
that $\M_p(d)$ is in any $L^r$, $r<p$ and worked with martingale differences
(not differences of reverse martingales).
But his argument applies to obtain the above proposition. We give the proof
of \eqref{weakMZ} in the \hyperref[appA]{Appendix}, for completeness. The proof of
\eqref{MZ}
is done in \cite{Woyczinski}. The argument is very similar
to the scalar case. Actually by the Banach principle (see Proposition~\ref{banprin}), using \eqref{weakMZ}, it is enough to show \eqref{MZ}
in the scalar case, see for instance the proof of Theorem~\ref{inemaxLIL}.
\end{rem}


Next, we obtain a similar result for $\M_2$, from which we derive the
compact LIL for stationary martingale differences (or reverse martingale
differences).

\begin{theo}\label{inemaxLIL}
Let $\X$ be a $(2,D)$-smooth separable Banach space, for some $D\ge1$.
For every $1\le p <2$, there exists a constant $C_p\ge1$, such that
for every
$d\in L^2(\Omega,\F_1,\X)$ (resp.
every $d\in L^2(\Omega,\F^0,\X)$) with $\E_0(d)=0$ (resp. $\E^{1}(d)=0$),
we have
%
\begin{equation}
\label{weakLIL} \bigl\|\M_2(d)\bigr\|_{p,\infty} \le C_p D\|d
\|_{2,\X}.
\end{equation}
In particular, $(d\circ\theta^n)_{n\ge0}$ satisfies the compact
{LIL}. Moreover, if $\theta$ is ergodic,
%
\begin{eqnarray}
\label{norm}\limsup_{n} \frac{|S_n(d)|_\X}{\sqrt{2nL(L(n))} }= \sup
_{x^*\in\X^*,|x^*|_{\X^*}\le1} \bigl\|x^*(d)\bigr\|_2 \le\|d\|_{2,\X} \qquad \as
\end{eqnarray}
and the cluster set of $ (\frac{S_n(d) }{\sqrt{2nL(L(n))} }
)_{n\ge1}$
is $\as$ a fixed compact set whose description is given in Appendix~\textup{\ref{cluster}}.
\end{theo}

\begin{rem}\label{rempisier}
 Of course, \eqref{weakLIL} is equivalent to the fact
that, for every $1\le p <2$,
there exists $\tilde C_p$, such that $\|\M_2(d)\|_{p} \le\tilde C_p D
\|d\|_{2,\X}$. This bound has been obtained in \cite{Pisier1}, Th\'eor\`eme 1,
for i.i.d. variables with values in a Banach space of type 2. Moreover,
it follows from Remarque 2 and the proposition page 208 of \cite
{Pisier1}, that if
every sequence of i.i.d. variables in $L^2(\Omega,\X)$ satisfy the
bounded LIL,
the space $\X$ must be of type $p$, for every $1<p<2$.
\end{rem}



Now, we deduce an almost sure invariance principle (ASIP) from Theorem~\ref{inemaxLIL}.
We first give the notations to specify what we mean by an ASIP, in the
Banach space setting.

Recall, that we denote by $\X^*$ the topological dual of $\X$. Let $X
\in L^2(\Omega,\X)$ such that \mbox{$\E(X)=0$}.
We define a bounded \emph{symmetric} bilinear operator $\K=\K_X$
from $\X^*\times\X^*$
to $\mathbb{R}$, by
\begin{eqnarray*}
\K\bigl(x^*,y^*\bigr) =\E\bigl(x^*(X)y^*(X)\bigr) \qquad \forall x^*,y^*\in\X^* .
\end{eqnarray*}
The operator $\K_X$ is called the \emph{covariance operator}
associated with $X$.

\begin{defn}
We say that a random variable $W\in L^2(\Omega, \X)$ is \emph
{Gaussian} if, for every $x^*\in\X^*$,
$x^*(W)$ has a normal distribution. We say that a random variable $X\in
L^2(\Omega,\X)$ is \emph{pregaussian},
if there exists a Gaussian variable $W\in L^2(\Omega,\X)$ with
the same covariance operator,
that is, such that $\K_X=\K_W$.
\end{defn}

\begin{defn}
We say that $(X_n)_{n\ge0}$ satisfies the almost sure invariance
principle ({ASIP}) if, without changing its distribution,
one can redefine the sequence $(X_n)_{n\ge0}$ on a new probability
space on which there exists a
sequence $(W_n)_{n\ge0}$ of centered i.i.d. Gaussian variables, such that
\[
\bigl|X_0+\cdots+X_{n-1}-(W_0+\cdots+W_{n-1})\bigr|_\X=\mathrm{o}
\bigl(\sqrt{nL\bigl(L(n)\bigr)}\bigr)\qquad  \as
\]
We shall say that $(X_n)_{n\ge0}$ satisfies the ASIP of covariance $\K
$, when $\K=\K_{W_0}$ is identified.
\end{defn}

We now recall an important result of Berger on the ASIP for martingale
differences.

\begin{prop}[(Berger {\cite{Berger}, Theorem~3.2})]\label{Berger}
Let $\X$ be a separable Banach space. Assume that $\theta$ is
ergodic. Let $d\in L^2(\Omega,\F_1,\X)$, with $\E_0(d)=0$. Assume
that $d$ is pregaussian and that $(d\circ\theta^n)_{n\ge0}$ satisfies
the {CLIL}. Then, for every $Y\in L^2(\Omega,\X)$, such that
$|S_n(Y)|_\X=\mathrm{o}(\sqrt{nL(L(n))})\ \as$, $((d+Y)\circ\theta
^n)_{n\ge0}$ satisfies the {ASIP} of covariance $\K_d$.
\end{prop}

Actually, Berger proved his result in the particular case where
$Y=Z-Z\circ\theta$ for some $Z\in L^2(\Omega,\X)$, but the proof applies
in the slightly more general situation above.

By \cite{LTbook}, Proposition~9.24, on any Banach space $\X$ of type
2 (in particular, on any
$2$-smooth Banach space), every $X\in L^2(\Omega, \X)$ is
pregaussian. Hence, Berger's result applies
as soon as the CLIL is satisfied and we deduce the following corollary.

\begin{cor}\label{corinemaxLIL}
Let $\X$ be a $2$-smooth separable Banach space. Assume that $\theta$
is ergodic.
For every $d\in L^2(\Omega,\F_1,\X)$, with $\E_0(d)=0$, $(d\circ
\theta^n)_{n\ge0}$ satisfies the {ASIP} of covariance $\K_d$.
\end{cor}

\begin{rem}\label{revrem}
Assume that $\dim \X=1$
and that $\theta$ is
ergodic. It follows from Corollary~2.5 of \cite{CM} that for $d\in
L^2(\Omega,\F^0,\X)$ such that $\E^1(d)=0$, $(d\circ\theta
^n)_{n\ge0}$ satisfies the ASIP. We do not know whether the ASIP holds
when $\dim \X\ge2$. The proof of
Proposition~\ref{Berger} given in \cite{Berger} does not seem to pass
to \emph{reverse} martingale differences.
\end{rem}
%

\subsection{Results for not necessarily adapted stationary
processes}\label{hannan}

We assume all along this subsection that $\theta$ is
invertible and
bi-measurable, in which case
we extend our filtration to $(\F_n)_{n\in\mathbb{Z}}$. Then, we write
$\F_{-\infty}:=\bigcap_{n\in\mathbb{Z}} \F_n$, $\F_\infty:=
\bigvee_{n\in\mathbb{Z}}\F_n$, and for every $n\in\overline{\mathbb{Z}}$,
$\E_n(\cdot)=\E(\cdot|\F_n)$ and $P_n:=\E_n-\E_{n-1}$. We say
that a random variable
$X\in L^1(\Omega,\X)$ is \emph{regular} if $\E_{-\infty}(X)=0$
and $X-\E_\infty(X)=0$.

\begin{theo}\label{theoHannanp}
Let $1<p<r\le2$ and
$D>0$.
Let $\X$ be a $(r,D)$-smooth separable Banach space and $X\in
L^p(\Omega,\F,\P,\X)$ be a regular variable. Assume moreover that
%
\begin{equation}
\label{Hannanp} \|X\|_{H_p}:=\sum_{n\in\mathbb{Z}}
\|P_n X\|_{p,\X}<\infty .
\end{equation}
Then, there exists (a universal) $C_{p,r}>0$, such that
%
\begin{equation}
\bigl\|\M_p(X)\bigr\|_{p,\infty} \le C_{p,r}D^{r/p} \| X
\|_{H_p} .
\end{equation}
Moreover,
%
\begin{equation}
\label{MZas} \bigl|S_n(X)\bigr|_\X/n^{1/p}\to0\qquad  \as
\end{equation}
\end{theo}

\begin{rem}
Theorem~\ref{theoHannanp} improves Corollary~1 of
\cite{Wu}, where \eqref{MZas} has been proved under a stronger condition
than \eqref{Hannanp}. The proof in \cite{Wu} is done for real-valued
variables but work in the above Banach setting as well.
\end{rem}

Now, we give a result under condition \eqref{Hannan}, which has been introduced
by Hannan \cite{Hannan}.

\begin{theo}\label{theoHannan}
Let $\X$ be a $(2,D)$-smooth separable Banach space, for some $D\ge1$.
Let $X\in L^2(\Omega,\F,\P,\X)$ be a regular random variable.
Assume moreover that
%
\begin{equation}
\label{Hannan} \|X\|_{H_2}:=\sum_{n\in\mathbb{Z}}
\|P_n X\|_{2,\X}<\infty .
\end{equation}
Then, for every $1\le p <2$, there exists (a universal) $C_p>0$, such that
%
\begin{equation}
\label{inemaxLILH} \bigl\|\M_2(X)\bigr\|_{p,\infty} \le C_pD \| X
\|_{H_2} .
\end{equation}
The series $d=\sum_{n\in\mathbb{Z}} P_1(X\circ\theta^n)$ converges
in $L^2(\Omega,\F_1,\X)$ and $\E_0(d)=0$. Moreover, writing $M_n:=
\sum_{k=0}^{n-1} d \circ\theta^k$,
we have
%
\begin{equation}
\label{martapp} |S_n-M_n|_\X= \mathrm{o}\bigl(\sqrt{nL
\bigl(L( n)\bigr)}\bigr)\qquad  \as
\end{equation}
%
\end{theo}

\begin{rem}
Theorem~\ref{theoHannan} improves Theorem~2 of
Wu \cite{Wu}, Theorem~2.1 of Liu and  Lin \cite{LL} (for $p=2$) and
Corollary~5.3 of Cuny \cite{Cuny1}. In \cite{Wu,LL,Cuny1} the authors prove
\eqref{martapp} under stronger conditions than \eqref{Hannan} and the proof
do not apply to infinite dimensional Banach spaces.
\end{rem}

In particular, we deduce the following corollary from Theorem~\ref{theoHannan},
Theorem~\ref{inemaxLIL} and Proposition~\ref{Berger}.

\begin{cor}\label{corberg}
Under the assumptions of Theorem~\ref{theoHannan}, $(X\circ\theta^n
)_{n\ge0}$ satisfies the CLIL
and the ASIP of covariance $\K_d$, where, for every $x^*\!,y^*\in\X
^*\!$, $\K_d(
x^*,y^*)$, $\K_d=\sum_{n\in\mathbb{Z}} \E(x^*(X_n)y^*(X))$.
Moreover, since, $\|d\|_{2,\X}\le\|X\|_{H_2}$,
\[
\limsup_n \frac{|S_n(X)|_\X}{\sqrt{2n L(L(n))}} \le\|X\|_{H_2}\qquad  \as
\]
\end{cor}

In order to check \eqref{Hannan} or \eqref{Hannanp}, it may be easier
to use the condition \eqref{hanbis} below.

\begin{lem}\label{Hanbis2}
Let $1< p\le2$. Let ${\mathcal H}$ be a separable real Hilbert space. Assume that
%
\begin{equation}
\label{hanbis} \sum_{n\ge1} \frac{\|\E_{-n}(X) \|_{p,\H}}{\sqrt n}<\infty
\quad\mbox{and} \quad \sum_{n\ge1} \frac{\|X-\E_{n}(X) \|_{p,\H}}{\sqrt n}<\infty .
\end{equation}
Then $X$ is regular and $\sum_{n\in\mathbb{Z}}\|P_n X\|_{p,\H
}<\infty$.
\end{lem}

\subsection{Functionals of Markov chains}
The situation considered in the previous paragraph includes the case of
stationary (ergodic) Markov chains. Let $Q$ be a
transition probability on a measurable space $(\SS,\S)$ admitting an
invariant probability $m$. Let
$(\Omega,\F,(\F_n)_{n\in\mathbb{Z}},\P, (W_n)_{n\in\mathbb
{Z}})$ be the canonical Markov chain associated with $Q$, that is,
$\Omega=\S^{\mathbb{Z}}$, $\F=\S^{\otimes\mathbb{Z}}$,
$(W_n)_{n\in\mathbb{Z}}$
the coordinates,
$\F_n=\sigma\{\ldots,W_{n-1},W_n\}$,
$\P\circ W_0^{-1}=m$ and $\P(W_{n+1}\in A |\F_n)=Q(W_n,A)$. Finally,
denote by $\theta$ the shift on $\Omega$.

Recall that $Q$ induces an operator on $L^2(\SS,m)$ that we still
denote by $Q$. If $\H$ is a separable real Hilbert space,
we denote by $\mathbf{Q}$ the analogous operator on $L^2(\SS,m,\H)$.
In particular,
for every $f\in L^2(\Omega,\H)$ and every $h\in\H$,
$\langle\mathbf{Q} f,h\rangle_\H= Q(\langle f ,h\rangle_\H)$.

Theorem~\ref{theoHannan} applies to that setting with $X=f(W_0)$,
where $f\in L^2(\SS,\H)$. Using Lem\-ma~\ref{Hanbis2}, it suffices to
check \eqref{hanbis}.
In that situation, the process is
adapted, that is, $X_0$ is $\F_0$-measurable. Hence, the second part
of condition \eqref{hanbis} is automatically satisfied while the first
part reads as follows
%
\begin{equation}
\label{Markov} \sum_{n\ge1}\frac{\|\mathbf{Q}^nf\|_{2,\H}}{\sqrt n} <\infty .
\end{equation}



%


\subsection{Results for non-invertible dynamical systems}

Here, we assume that $\theta$ is non-invertible. Let us write $\F
^n=\theta^{-n}(\F)$, for every $n\ge0$. Denote $\F^\infty=\bigcap
_{n\ge0}\F^n$.

In this case, there exists a Markov operator $K$, known as the
Perron--Frobenius operator, defined by
%
\begin{equation}
\label{def} \int_\Omega X (Y\circ\theta) \,\mathrm{d}\P=\int
_\Omega(KX) Y \,\mathrm{d}\P\qquad  \forall X,Y\in L^2(\Omega, \F,\P)
.
\end{equation}

Then, we have for every $X\in L^1(\Omega,\F^0,\P)$,
%
\begin{equation}
\label{relfro} \E^n(X)=\bigl(K^nX\bigr)\circ
\theta^n .
\end{equation}
If $\H$ is a separable real Hilbert space, we extend $K$
to $L^2(\Omega,\F,\P,\H)$, in a way similar to \eqref{def}. We
denote by $\mathbf{K}$ the obtained operator.

\begin{theo}\label{dynsys}
Let $(\Omega,\F,\P,\theta)$ be a non-invertible dynamical system.
Let $X\in L^2(\Omega,\H)$ be such that
%
\begin{equation}
\label{conddynsys} \sum_{n\ge0} \frac{\|\mathbf{K}^nX\|_{2,\H}}{\sqrt n} <\infty .
\end{equation}
Then, for every $1<p<2$, there exists $C_p>0$ such that
\[
\bigl\|\M_{2}(X)\bigr\|_{p,\H}\le C_p \sum
_{n\ge0} \frac{\|\mathbf{K}^nX\|
_{2,\H}}{\sqrt n} .
\]
Moreover, there exists $d\in L^2(\Omega,\F,\P,\H)$ with $\E
^1(d)=0$, such that,
writing $M_n:=\sum_{k=0}^{n-1}d\circ\theta^k$, we have
%
\begin{equation}
\label{martappdynsys} |S_n-M_n|_{\H}=\mathrm{o}\bigl(\sqrt{nL
\bigl(L(n)\bigr)}\bigr)\qquad  \as
\end{equation}
\end{theo}

\begin{rem}
 It follows from \eqref{martappdynsys}
that $(X\circ\theta^n)_{n\ge0}$ satisfies the {CLIL}, but we do not
know whether it satisfies the {ASIP} in general, except when $\H$ has
dimension one (see Remark~\ref{revrem}).
\end{rem}

\section{Applications, examples}\label{examples}

Now, we give several applications of the previous results. We do not intend
to give all possible examples where our conditions apply, but we try to
provide examples illustrating the different situations we have considered.

For instance, our results on the Marcinkiewicz--Zygmund strong laws
(and on the LIL) may be used (in the one-dimensional case) to
obtain almost-sure invariance principles with rate as in \cite{Wu}
(see also
\cite{DDM} or \cite{CM}).

We start with a one-dimensional situation.

\subsection{\texorpdfstring{$\phi$}{phi}-mixing sequences}

Let us recall the definition of the $\phi$-mixing coefficients,
introduced by
Dedecker and Prieur \cite{DP}. Examples of $\phi$-mixing sequences
may be found there as well.

\begin{defn}\label{defphi}
For any integrable random variable $X$, let us write
$X^{(0)}=X- \E(X)$.
For any random variable $Y$ with values in
${\mathbb R}$ and any $\sigma$-algebra $\F$, let
\[
\phi(\F, Y)= \sup_{x \in{\mathbb R}} \bigl\| \E \bigl( (\I_{Y \leq x})^{(0)}
| \F \bigr)^{(0)} \bigr\|_\infty.
\]
For a sequence $\mathbf{Y}=(Y_i)_{i \in{\mathbb Z}}$, where $Y_i=Y
\circ\theta^i$ and $Y$ is an $\F_0$-measurable and real-valued
random variable, let
\begin{eqnarray*}
\phi_{ \mathbf{Y}}(n) = \sup_{ i\ge n} \phi(\F_0,
Y_{i}) .
\end{eqnarray*}
\end{defn}

We need also the following technical definition.

\begin{defn}
If $\mu$ is a probability measure on $\mathbb R$ and $p \in\,]1,
\infty)$, $M \in(0, \infty)$, let $\Mon_p(M,\mu)$ denote the set of
functions $f\dvtx \R\to\R$ which are monotone on some interval and null
elsewhere and such that $\mu(|f|^p)\leq M^p$. Let $\Mon^c_p(M,\mu)$ be
the closure in ${\mathbb L}^p(\mu)$ of the set of functions which can
be written as $\sum_{\ell=1}^L a_\ell f_\ell$, where $\sum_{\ell
=1}^L|a_\ell| \leq1$ and $f_\ell\in\Mon_p(M, \mu)$.
\end{defn}

\begin{theo}\label{asipgenephi}
Let $X = f(Y) - \E( f(Y))$, where $Y$ is an $\F_0$-measurable random
variable. Let $P_{Y}$ be the distribution of $Y$ and $p \in\,]1,\infty
]$. Assume that $f$ belongs to $\Mon^c_p(M,P_{Y})$ for some $M>0$, if
$2\le p <\infty$ and that
$f$ has bounded variation if $p=\infty$. Assume moreover that
%
\begin{equation}
\label{condDDM} \sum_{k \geq1} \frac{ \phi^{(p-1)/p}_{ \mathbf{Y}}(k)}{k^{1/2}} <
\infty.
\end{equation}
Then, if $1<p<2$, $(X\circ\theta^n)_{n\in\mathbb{Z}}$ satisfies the
conclusion of Theorem~\ref{theoHannanp} and if $p\ge2$, $(X\circ\theta^n)_{n\in\mathbb
{Z}}$ satisfies the conclusion of Theorem~\ref{theoHannan}.
\end{theo}

\begin{rem}
 When $p=2$, Dedecker, Gou\"ezel and Merlev\`ede
\cite{DGM} proved that the\break condition $\sum_{k \geq1}
k^{1/\sqrt3-1/2} \phi^{1/2}_{ \mathbf{Y}}(k) < \infty$ implies that
$\sum_{n\ge1} \P(\max_{1\le k \le2^n}
|S_k|>\break  C 2^{n/2}(L(n))^{1/2})<\infty$ (which implies the bounded LIL).
\end{rem}

\begin{pf*}{Proof of Theorem~\ref{asipgenephi}} Assume first that $1< p<\infty$. Since $f\in\Mon
_p^c(M, P_{Y_0})$, there exists
a sequence of functions
\[
f_L=\sum_{k=1}^L
a_{k,L} f_{k,L} ,
\]
such that for every $L\ge1$, $\sum_{k=1}^L |a_{k,L}|\le1$, for every
$1\le k \le L$, $f_{k,L}$ is monotonic on some interval and
null elsewhere, and $\|f_{k,L}(Y_0)\|_p\le M$ and such that
$(f_L)_{L\ge1}$
converges in $L^p(P_{Y_0})$ to $f$.
Hence,
\begin{eqnarray*}
&&\bigl\|\E_0\bigl(f(Y_n)\bigr)-\E\bigl(f(Y_n)
\bigr)\bigr\|_p \\
&&\quad =\lim_{L\to\infty} \bigl\|\E_0
\bigl(f_L(Y_n)\bigr)-\E\bigl(f_L(Y_n)
\bigr)\bigr\|_p
\\
&&\quad \le\liminf_{L\to\infty}\sum_{k=1}^L|a_{k,L}|
\bigl\|\E_0\bigl(f_{k,L}(Y_n)\bigr)- \E
\bigl(f_{k,L}(Y_n)\bigr)\bigr\|_p \le C_p
M \phi^{(p-1)/p}_{ \mathbf{Y}}(n) ,
\end{eqnarray*}
where we used Lemma~5.2 of \cite{DGM1} for the last estimate.

To conclude in that case, we notice first that we are in the adapted case,
and that Theorem~\ref{theoHannanp} (when $1<p<2$) and Theorem~\ref
{theoHannan} (when $p\ge2$) apply by Lemma~\ref{Hanbis2}.

Assume that $p=\infty$ and that $f$ has bounded variation.
Hence $f$ is the difference of two monotonic functions, to which
we apply Lemma~5.2 of \cite{DGM1} with $p=\infty$. Then, we conclude
as above.
\end{pf*}

\subsection{\texorpdfstring{$\X$}{X}-valued linear processes}

Let $(\Omega,\F,\P)$ be a probability space and $\theta$ be an
ergodic invertible and bi-measurable transformation on $\Omega$. Let
$\X$ be a separable $r$-smooth Banach space, for some $1<r\le2$. Let
$\xi\in L^p(\Omega,\F_0,\P, \X)$ for some $p>1$. Assume that
$\E(\xi|\F_{-1})=0$ and define
$\xi_n=\xi\circ\theta^n$, $n\in\mathbb{Z}$.


Let $(A^{(k)})_{k\in\mathbb{Z}}$ be a (not necessarily stationary)
sequence of
random variables with values in $L^\infty(\Omega,\F_{k-1},\B(\X
))$, where
$\B(\X)$ stands for the Banach space of bounded
(linear) operators on $\X$. For every $k,n\in\mathbb{Z}$, define
$A^{(k)}_n=A^{(k)}\circ\theta^n $.
Assume that
%
\begin{eqnarray}
\label{hanlin} \sum_{k\in\mathbb{Z}} \bigl\|A^{(k)}
\bigr\|_{\infty, \B(\X)} <\infty .
\end{eqnarray}
Then, the process
\[
X_n:= \sum_{k\in\mathbb{Z}}A_n^{(k)}
\xi_{n+k} , \qquad n\in\mathbb{Z}
\]
is well defined in $L^p(\Omega,\X)$ and is stationary.

\begin{cor}\label{linear}
Assume that $1<p<r\le2$ or $p=r=2$. Let $(X_n)$ be a linear process as above.
Then,
%
\begin{equation}
\sum_{n\in\mathbb{Z}} \|P_nX_0
\|_{p,\X} <\infty .
\end{equation}
Hence, Theorem~\ref{theoHannanp} applies when $1<p<2$ and Theorem~\ref{theoHannan} applies when $p=2$.
\end{cor}

\subsection{Functions of real-valued linear processes}

Let $(\xi_n)_{n\in\mathbb{Z}}$ be a sequence of \emph{independent}
identically distributed \emph{real}
random variables in $L^2(\Omega,\F,\P)$. Let $(a_n)_{n\in\mathbb{Z}}$
be in $\ell^1$. We consider a linear process defined by
\[
Y_n := \sum_{k\in\mathbb{Z}} a_k
\xi_{n-k} \qquad \forall n\in\mathbb{Z} .
\]
For every $n\in\Z$, write $\F_n=\sigma\{\ldots, \xi_{n-1},\xi
_{n}\}$.

We denote by $\Lambda$ the class of non-decreasing continuous and bounded
functions on
$[0,+\infty[$, such that $\varphi(0)=0$, and satisfying one of the following
\begin{eqnarray*}
&\mbox{$\varphi^2$ is concave};&
\\
&\varphi(x)=C\min\bigl(1,x^\alpha\bigr) \qquad \forall x\ge0 , \mbox{ for some
$0< \alpha\le1$, $C>0$}.&
\end{eqnarray*}

Let $r\ge1$. Let $f$ be
a real valued function such that
%
\begin{equation}
\label{condcont0} \bigl|f(x)-f(y)\bigr|\le\varphi\bigl(|x-y|\bigr) \bigl(1+|x|^r+|y|^r
\bigr) \qquad \forall x,y\in\R .
\end{equation}
Our functions are unbounded and their continuity
is locally controlled by $\varphi$.

We want to study the process $(X_n)_{n\in\Z}$ given by
\[
X_n:=f(Y_n)-\E\bigl(f(Y_n)\bigr)\qquad  \forall n
\in\Z .
\]

\begin{cor}\label{flp}
Let $\varphi\in\Lambda$ and $r\ge1$. Let $\xi_0\in L^{2r}(\Omega
,\F,\P)$ and $f$ satisfy \eqref{condcont0}.
Let $(a_n)_{n\in\mathbb{Z}}\in\ell^1$. Consider the process
$(X_n)_{n\ge0}$ above.
If
\[
\sum_{n\ge1} \varphi\bigl(|a_n|\bigr) <\infty \quad \mbox{or}\quad
\sum_{n\ge1} \frac{\varphi( \sum_{k\ge n}|a_k|)}{\sqrt n} <\infty ,
\]
then $(X_n)_{n\ge0}$ satisfies the conclusion of Theorem~\ref{theoHannan}.
\end{cor}

We give the proof in the \hyperref[appA]{Appendix}.

\begin{rem}
Notice that condition (3.1) of \cite{LL} implies
\eqref{condcont0} with $\varphi(x)=\min(1,x)$. Hence, Corollary~\ref{flp}
improves Corollary~3.1 of \cite{LL} when $p=2$.
\end{rem}





\subsection{A non-adapted example}

We now consider an example of a non-adapted process for which new ASIP
with rates have been obtained very recently, see
Dedecker, Merlev\`ede and P\`ene \cite{DMP} and the references therein.

Let $d\ge2$ and $\theta$ be an ergodic automorphism of the
$d$-dimensional torus $\Omega=\Omega_d=\mathbb{R}^d/\mathbb{Z}^d$. Denote
by $\F$ the Borel $\sigma$-algebra of $\Omega$
and take $\P$ to be the Lebesgue measure on $\Omega$.

For every $\mathbf{k}=(k_1,\ldots, k_d)\in\mathbb{Z}^d$, write
$|\mathbf{k}|:=
\max_{1\le i\le d}|k_i|$. If $\H$ is a Hilbert space and if $f\in L^2
(\Omega,\F,\P,\H)$, we denote by $(c_\mathbf{k})_{\mathbf{k}\in
\mathbb{Z}^d}
=(c_{\mathbf{k},\H})_{\mathbf{k}\in\mathbb{Z}^d}$ its Fourier coefficients,
that is, $c_{\mathbf{k},\H}= \int_{[0,1]^d} f(x) \mathrm{e}^{-2\mathrm{i}\pi
\langle x,\mathbf{k}
\rangle_d}\P( \mathrm{d}x)$, for every $\mathbf{k}\in\mathbb{Z}^d$, where
$\langle\cdot,
\cdot\rangle_d$ stands for the inner product on $\mathbb{R}^d$.

\begin{cor}\label{DMP}
Let $\H$ be a Hilbert space and $f\in L^2(\Omega, \H)$. Assume that
there exists $\beta>2$ and $C>0$ such that
\[
\sum_{|\mathbf{k}|\ge m} |c_\mathbf{k}|_\H^2
\le\frac{C}{L(
m)(L(L( m)))^\beta}\qquad  \forall m\ge1.
\]
Then, $(f\circ\theta^n)_{n\ge0}$ satisfies the ASIP with covariance
operator given by
$\K(x,y):=\linebreak[4]  \sum_{m\in\mathbb{Z}} \E(\langle x,f\rangle_\H\langle
y,f\circ\theta^n\rangle_\H)$,
for every $x,y\in\H$.
\end{cor}

\begin{rem}
 Dedecker, Merlev\`ede and P\`ene
\cite{DMP}, Theorem~2.1, obtained the ASIP when $\H=\mathbf{R}^m$ and their
condition
requires $\beta>4$. When $m=1$, rates in the ASIP are also provided in
\cite{DMP}.
\end{rem}

\begin{pf*}{Proof of Corollary~\ref{DMP}} It follows from the proof of Propositions 4.2 and
4.3 of \cite{DMP} (notice that the proofs work in the Hilbert space setting)
that there exists a filtration $(\F_n)_{n\in\mathbb{Z}}$ (defined at
the beginning of paragraph 3 of \cite{DMP}) such that
$\F_{n}=\theta^{-n}(\F_0)$ and
\[
\bigl\|\E_{-n}(f)\bigr\|_{2,\H} =\mathrm{O}\biggl(\frac{1}{\sqrt{nL(n)^\beta}}\biggr)
\quad \mbox{and}\quad  \bigl\|\E_n(f)-f\bigr\|_{2,\H} =\mathrm{O}\biggl(\frac{1}{\sqrt{nL(n)^\beta}}
\biggr) .
\]
Then, the result follows from Lemma~\ref{Hanbis2}. \end{pf*}

\subsection{Cramer--von Mises statistics}

We use our previous notations, see Section~\ref{hannan}.

Let $Y\in L^0(\Omega,\F_0,\P)$. For every $n\in\mathbb{Z}$, let
$Y_n:= Y\circ\theta^n$ and $X_n:= t\mapsto\I_{Y_n\le t}-F(t)$,
where $F(t)=\P(Y\le t)$.

Let $1<r\le2$. For every $\sigma$-finite Borel measure $\mu$ on
$\mathbb{R}$, we may see
$(X_n)_{n\in\mathbb{Z}}$ as a process with values in the $r$-smooth
Banach space
$L^r(\mathbb{R},\mu)$, as soon as
%
\begin{equation}
\label{condminCM} \int_0^\infty\bigl(1-F(t)
\bigr)^r \mu(\mathrm{d}t) + \int_{-\infty}^0
F(t)^r\mu(\mathrm{d}t) <\infty ,
\end{equation}
which is satisfied whenever $\mu$ is finite.

Define $F_\mu$ by $F_\mu(x)=-\mu([x,0[)$ if $x\le0$ and $F_\mu(x)=
\mu([0,x[)$ if $x\ge0$. Let $1< p\le2$. Then, under \eqref{condminCM},
$X_0\in L^p(\Omega,L^r(\mu))$
if and only if
%
\begin{equation}
\label{condminCMbis} \E\bigl(\bigl|F_\mu(Y_0)\bigr|^{p/r}\bigr)<
\infty .
\end{equation}

We want to understand the asymptotic behaviour of the process
$F_n=S_n(X)/n$ (with values in $L^2(\mathbb{R},\mu)$), and more
particularly of $D_n(\mu):=\|F_n\|_{2,\mu}$.
When $\mu=P_Y=\P\circ Y^{-1}$, $D_n(\mu)^2$ is known as the
Cramer--von Mises statistics.

It follows from Lemma~\ref{Hanbis2}, that if
$(X_n)_{n\in\mathbb{Z}}$ satisfies
%
\begin{equation}
\label{condCM} \sum_{n\ge1} \frac{ (\E(\| \E_{-n}(X_0)\|_{2,\mu}^p) )^{1/p}} {
n^{1/2}} <\infty
,
\end{equation}
for some $1<p\le2$, then $(X_n)_{n\in\mathbb{Z}}$ satisfies
Theorem~\ref{theoHannanp} if $1<p<2$ and Theorem~\ref{theoHannan} if $p=2$.
Hence, we have
the following corollary.

\begin{cor}\label{corCM}
Let $1<p<r\le2$ or $p=r=2$. With the above notations, assume that
\eqref{condminCM},
\eqref{condminCMbis} and
\eqref{condCM} be satisfied. Then,
\begin{eqnarray*}
\lim_n n^{1-1/p} D_n(\mu) &=& 0\qquad  \as
\qquad \mbox{if $1<p<2$} ;
\\
\limsup_n \frac{n^{1/2}}{(2L(L(n)))^{1/2}} D_n(\mu)& =&
\Lambda_\mu\qquad  \as \qquad \mbox{if $p=2$} ,
\end{eqnarray*}
where $\Lambda_\mu^2:={ \sup_{\|f\|_{2,\mu,\mathbb{R}}\le1}}
\int_{\mathbb{R}^2} f(s)f(t) C(s,t) \mu(\mathrm{d}s)\mu(\mathrm{d}t)$ and $C(s,t):=
{ \sum_{n\in\mathbb{Z}}} ( \P(Y_0\le s,\allowbreak  Y_n\le t)
-F(s)F(t) )$.
\end{cor}

\begin{pf} Apply Lemma~\ref{Hanbis2}, Theorem~\ref{theoHannanp} and
Theorem~\ref{theoHannan}. The expression of $\Lambda_\mu^2$ follows, for
instance, from Proposition~1 of Merlev\`ede \cite{Merlevede}.
\end{pf}

In the context of $\phi$-mixing sequences, when $\mu$ is finite, Corollary~\ref{corCM} applies as soon as
$\sum_{n\ge1} \frac{\phi_\mathbf{Y}(n)^{1/2}}{n^{1/2}}<\infty$.

Other examples where \eqref{condCM} is satisfied may be found in \cite{DMemp}.


\section{Proof of the results for Banach-valued martingales}\label{proofs1}

%


%


\begin{pf*}{Proof of Theorem~\ref{inemaxLIL}} Let us prove \eqref{weakLIL}.
We start with the case $d\in L^2(\Omega,\F_1,\P)$ and $\E_0(d)=0$.

When $d\in L^2(\Omega,\F^{0},\P)$ and $\E^1(d)=0$, the proof is the same,
with the obvious changes, noticing that for every $n\ge1$,
$(S_n(d)-S_{n-k}(d))_{0\le k \le n}$ is a $(\F^{n-k})_{0\le k \le
n}$-martingale
and
that $\max_{1\le k\le n}|S_k(d)|_\X\le2
\max_{1\le k\le n}|S_n(d)-S_{n-k}(d)|_\X$.


Clearly, by homogeneity, it suffices to prove the
result when $\|d \|_{2,\X}=1$. Let $\lambda>0$ and $1\le p <2$. Let
us prove that there exists
$C_p\ge1$, independent of $\lambda$ such that
%
\begin{equation}
\label{inmax} \lambda^p \P\bigl(M^* > \lambda\bigr) \le
D^p C_p ^p ,
\end{equation}
where
\begin{eqnarray*}
M^*=M^*(d):=\sup_{s\ge0}\frac{\max_{1\le k\le2^s} |S_k(d)|_\X} {
2^{s/2}(L( s))^{1/2}} .
\end{eqnarray*}
Since $\M_2(d)\le C M^*$, this will imply the desired result. Notice that
\eqref{inmax} holds trivially when $0<\lambda<D$. Assume then that
$\lambda\ge D$.

Let $S\ge1$ be an integer, fixed for the moment. For simplicity, we write
$S_n:=S_n(d)$.

We have, using Doob's maximal inequality for the submartingale
$(|S_n|_\X)_{n\ge1}$, and \eqref{smooth}
%
\begin{eqnarray}\label{firstest}
\P\biggl(\sup_{1\le s\le S}\frac{\max_{1\le k\le2^s} |S_k|_X} {
2^{s/2}(L( s))^{1/2}} > \lambda
\biggr)&\le& \frac{1}{\lambda^2}\sum_{s=1}^S
\frac{\E(\max_{1\le k\le2^s}
|S_k|_\X^2)} {
2^sL(s)}\nonumber
\\[-8pt]\\[-8pt]
&\le& \frac{2}{\lambda^2}\sum_{s=1}^S
\frac{ \E
(|S_{2^s}|_\X^2)} {
2^sL(s)} \le\frac{2 D^2 S}{\lambda^2 } .\nonumber
\end{eqnarray}


We make use of truncations. Let $\alpha>0$ be fixed for the moment.
Let us write $d_n:=d\circ\theta^{n-1}$, $n\ge1$.
For every $s\ge1$, $k\ge1$ define
\begin{eqnarray*}
e_k^{(s)}&:=& d_k\mathbf{1}_{\{|d_k |_\X\le\alpha\lambda2^{s/2}
/(L(s))^{1/2}\}} ;\qquad
d_k^{(s)}:= e_k^{(s)} -\E
\bigl(e_{k}^{(s)}|\F_{k-1}\bigr) ;\qquad  \tilde
d_k^{(s)}:=d_k-d_k^{(s)},
\\
S_k^{(s)}&:=&\sum_{i=1}^k
d_i^{(s)} ; \qquad \tilde S_k^{(s)}:=
S_k-S_k^{(s)},
\\
T_s&:=&4 \sum_{i=1}^{2^s} \E
\bigl(|d_i|_\X^2|\F_{i-1}\bigr) ;\qquad
T_s^{(s)}:= \sum_{i=1}^{2^s}
\E\bigl(\bigl|d_i^{(s)}\bigr|_\X^2|
\F_{i-1}\bigr) .
\end{eqnarray*}
Notice that, for every $s\ge1$,
%
\begin{equation}
\label{condvar} T_s^{(s)}\le T_s .
\end{equation}
Let $\beta>0$ be fixed for the moment. Define the events
\begin{eqnarray*}
A_s&:=&\biggl\{\frac{\max_{1\le k\le2^s} |S_k|_\X} {
2^{s/2}(L(s))^{1/2}} > \lambda\biggr\} ;\qquad
B_s:= \biggl\{\frac{\max_{1\le k\le2^s} |S_k^{(s)}|} {
2^{s/2}(L( s))^{1/2}} > \lambda/2\biggr\},
\\
C_s&:=& \biggl\{\frac{\max_{1\le k\le2^s} |\tilde S_k^{(s)}|_\X} {
2^{s/2}(L( s))^{1/2}} > \lambda/2\biggr\} ;\qquad
D_s:=\biggl\{\frac{T_s}{2^s}>\beta\lambda^2\biggr\} ;\qquad
E_s:=B_s \cap\biggl\{\frac{T_s^{(s)}}{2^s}\le\beta
\lambda^2\biggr\} .
\end{eqnarray*}
Using \eqref{condvar}, we see that $B_s\cap D_s^c\subset E_s$.
In particular, we have
\begin{eqnarray*}
A_s\subset B_s \cup C_s ;\qquad
B_s\subset D_s \cup E_s .
\end{eqnarray*}
Hence,
\begin{eqnarray*}
\biggl\{\sup_{s\ge S}\frac{\max_{1\le k\le2^s} |S_k|_\X} {
2^{s/2}(L( s))^{1/2}} > \lambda\biggr\} =
\bigcup_{s\ge S}A_s \subset \biggl(\bigcup
_{s\ge S} C_s\biggr) \cup
\biggl(\bigcup_{s\ge S} D_s\biggr)\cup
\biggl(\bigcup_{s\ge S}E_s\biggr)
.
\end{eqnarray*}
Now, $\bigcup_{s\ge S}D_s=\{\sup_{s\ge S}\frac{T_s}{2^s}>\beta
\lambda^2\}$,
hence by Hopf maximal inequality \eqref{Hopf},
using that $\E(|d_1|_\X^2)=1$,
%
\begin{equation}
\label{hopf} \P\biggl(\bigcup_{s\ge S}D_s
\biggr)\le\P\biggl(\bigcup_{s\ge1}D_s
\biggr)\le \frac{4}{\beta\lambda^2} .
\end{equation}

We also easily see that, interverting $\sum$ and $\E$ in
\eqref{borel},
%
\begin{eqnarray}\label{borel}
\P\biggl(\bigcup_{s\ge S} C_s
\biggr)&\le&\frac{2}{\lambda} \sum_{s\ge0}
\frac{\E(\max_{1\le k\le2^s} |\tilde S_k^{(s)}|_\X)} {
2^{s/2}(L( s))^{1/2}}\nonumber
\\[-8pt]\\[-8pt]
&\le&\frac{4}{\lambda} \sum_{s\ge1}
\frac{2^{s/2}}{(L(s))^{1/2}} \E\bigl(|d_1|_\X \mathbf{1}_{\{ |d_1|_\X\ge\alpha\lambda2^{s/2}/(L(s))^{1/2}\}}\bigr)
\le \frac{4C}{\alpha\lambda^2} ,\nonumber
\end{eqnarray}
where we also used that there exists $C>0$ such that for every $u>0$,
\begin{eqnarray*}
\sum_{ s \le u} 2^{s/2}/\bigl(L(s)
\bigr)^{1/2} \le\sum_{ s\le\sqrt u} 2^{s/2}
+ \frac{1}{(L(\sqrt u))^{1/2}}\sum_{\sqrt u<s \le u} 2^{s/2} \le
C 2^{u/2}/L(u)^{1/2} .
\end{eqnarray*}

It remains to deal with $\bigcup_{s\ge S} E_s$. We need the following lemma
from Dedecker, Gou\" ezel and Merlev\`ede \cite{DGM}, Proposition A.1
(see also
Merlev\`ede \cite{Merlevede}, Lemma~1), whose proof follows from
Pinelis \cite{Pinelis}, Theorem~3.4. The proof in \cite{DGM}
is done in the scalar case (and in \cite{Merlevede} in the Hilbert case)
but it easily extends to $2$-smooth Banach spaces, since Theorem~3.4 in
\cite{Pinelis} is proved in that
setting. A related inequality in the scalar case is stated in Freedman
\cite{Freedman}, Theorem~1.6.

\begin{lem}
Let $\X$ be a $(2,D)$-smooth Banach space. Let $c>0$. Let $(\F
_j)_{j\ge0}$ be a non-decreasing filtration and
$(d_j)_{j\ge1}$ a sequence of random variables adapted to $(\F
_j)_{j\ge0}$,
such that for every $j\ge1$, $|d_j|_\X\le c$ a.s. and $\E(d_j|\F_{j-1})=0$
a.s. Then, for all $x,y>0$ and all integer $n\ge1$, we have
%
\begin{equation}
\label{expin} \P\Biggl( \max_{1\le k\le n} \Biggl|\sum
_{i=1}^kd_i\Biggr|_\X>x ; \sum
_{i=1}^n \E\bigl(|d_i|_\X^2|
\F_{i-1}\bigr) \le y/D^2 \Biggr) \le2 \exp \biggl(-
\frac{y}{c^2}h\biggl(\frac{xc}{y}\biggr) \biggr) ,
\end{equation}
where $h(u)=(1+u)\log(1+u)-u$.
\end{lem}

Let $s\ge S$. Let us apply the lemma to the sequence of martingale differences
$(d_i^{(s)})$ (in this case, we may take $c=2\alpha\lambda2^{s/2}/
(L(s))^{1/2}$), with $x=\lambda2^{s/2-1}(L(s))^{1/2}$,
$y= \beta D^2 \lambda^2 2^s$ and $n=2^s$. We obtain, taking $\alpha
=D^2\beta$,
\begin{eqnarray*}
\P(E_s)\le2 \exp \biggl( - \frac{D^2 \beta L(s)}{4\alpha^2}h\biggl(
\frac
{\alpha} {
D^2\beta}\biggr) \biggr)= 2 \exp \biggl( - \frac{ L(s)h(1)}{4D^2\beta} \biggr) =
\frac{2}{s^{h(1)/4D^2\beta}} .
\end{eqnarray*}
Hence, if $h(1)/(4D^2\beta) >1$, we see that
\begin{eqnarray*}
\sum_{s\ge S} \P(E_s)\le
\frac{2}{(h(1)/4D^2\beta
-1)S^{h(1)/4D^2\beta-1}} .
\end{eqnarray*}
Take $\beta= \frac{(2-p)h(1)}{8D^2}$ and $S=[\lambda^{2-p}]$. Then,
$h(1)/4D^2\beta-1= 2/(2-p)-1=p/(2-p)$ and
%
\begin{equation}
\label{lastest} \sum_{s\ge S} \P(E_s) \le
\frac{C}{(2-p)\lambda^p} .
\end{equation}
Recall that we
assume that $\lambda\ge D$, in particular $\frac{1}{\lambda^2}\le
\frac{D^{p-2}}{\lambda^p}$. Combining \eqref{firstest},
\eqref{hopf}, \eqref{borel} and \eqref{lastest}, we infer that,
there exists $C>0$, such that
\[
\lambda^p\P\bigl(M^*>\lambda\bigr)\le\frac{C D^p }{2-p} ,
\]
which ends the proof of \eqref{weakLIL}.

Let us prove that $(d\circ\theta^n)_{n\in\mathbb{N}}$ satisfies the
{CLIL}. We shall use the Banach
principle, see Proposition~\ref{banprin}. By definition of the Bochner
spaces, there exists $(d^{(m)})_{m\ge1}$, converging in
$L^2(\Omega,\X)$ to $d$, such that for every $m\ge1$, there exist
$k_m\ge1$, $\alpha_1, \ldots, \alpha_{k_m}\in\X$
and $A_1,\ldots, A_{k_m}\in\F_1$ such that
\[
d^{(m)}=\sum_{i=1}^{k_m}
\alpha_i \mathbf{1}_{A_i}.
\]
Write $\tilde d^{(m)}:= d^{(m)}-\E_0(d^{(m)})$. Then, $(\tilde
d^{(m)})_{m\ge1}$ converges
in $L^2(\Omega,\X)$ to $d$. Hence, by the Banach principle, it
suffices to
prove that $(\tilde d^{(m)}\circ\theta^n )_{n\in\mathbb{N}}$
satisfies the CLIL
for every $m\ge1$. But, by construction, $(\tilde d^{(m)}\circ\theta
^n )_{n\in\mathbb{N}}$ is a stationary sequence of martingale
differences taking values in a \emph{finite} dimensional Banach space
(i.e., $\Vect
\{\alpha_i \dvt  1\le i\le k_m\}$), in which case the compact LIL and the
bounded LIL are equivalent. But the bounded LIL in that case follows from
\eqref{weakLIL}, hence the result.


It remains to prove \eqref{norm}. By the bounded LIL the variable
$\limsup_n
\frac{|S_n(d)|_\X}{\sqrt{2nL(L(n))}}$ is well-defined $\P$-a.s. and
must be $\theta$-invariant. By ergodicity, there exists
$S\ge0$, such that $\limsup_n \frac{|S_n(d)|_\X}{\sqrt {2nL(L(n))}}=S$\ \as\ Let $M:=\sup_{|x^*|_{\X^*}\le1}\|x^*(d)\|_2$.
Let us prove that $S=M$. Let $\varepsilon>0$. There exists
$x^*_\varepsilon\in\X^*$, with $|x^*_\varepsilon|_{\X^*}\le1$,
such that $\|x^*_\varepsilon(d)\|_2
\ge M-\varepsilon$. Since,
$|S_n(d)|_\X\ge|x^*_\varepsilon(S_n(d))|$, it follows from the LIL
for real-valued martingales (with stationary ergodic increments) that
\[
S\ge M-\varepsilon .
\]
Letting $\varepsilon\to0$, we see that $S\ge M$. Let us prove the
converse inequality.

Let $x^*\in\X^*$. By the LIL for real-valued, stationary and ergodic
martingale differences, $\limsup_n S_n(x^*(d))/\sqrt{2nL(L(n))}= \|
x^*(d)\|_2\ \as$
Hence, by the compact LIL and Proposition~\ref{Kuelbs},
there exists a compact set $K\in\X$, such that for $\P$-a.e. $\omega
\in\Omega$,
the cluster set of $\{S_n(d)(\omega)/\sqrt{2nL(L(n))}, n\ge1\}$ is
$K$. Let $x\in K$ be such that $|x|_\X=S$, and let
$x^*\in\X^*$ be such that $|x^*|_{\X^*}=1$ and
$x^*(x)=|x|_\X$. For $\P$-a.e. $\omega\in\Omega$,
there exists $(n_k=n_k(\omega))_{k\ge1}$ such that
$S_{n_k}(d)(\omega)\sqrt{2n_kL(L(n_k))}\displaystyle \mathop{\longrightarrow}_{k\to\infty}^{|\cdot|_\X} x$. In particular
\[
x^*\bigl(S_{n_k}(d) (\omega)\sqrt{2n_kL
\bigl(L(n_k)\bigr)}\bigr)\displaystyle \mathop{\longrightarrow}_{k\to\infty}^{|\cdot|_\X}
x^*(x)=S \le \limsup_n S_{n}\bigl(x^*(d)\bigr) (
\omega)\sqrt{2nL\bigl(L(n)\bigr)} .
\]
But, by the real LIL, for $\P$-a.e. $\omega\in\Omega$,
\[
\limsup_n S_{n}\bigl(x^*(d)\bigr) (\omega)
\sqrt{2nL\bigl(L(n)\bigr)} \le\bigl\|x^*(d)\bigr\|_2 \le M ,
\]
which ends the proof. \end{pf*}


%



\section{Proof of the results for stationary processes}\label{proofs2}

\subsection{Proof of Theorem \texorpdfstring{\protect\ref{theoHannan}}{2.10}}

Recall that we assume here $\theta$ to be invertible. Let $\X$ be a
$2$-smooth Banach space.

Define
%
\begin{equation}
\label{banhan} H_2:=\biggl\{ Z\in L^2(\Omega, \X) \dvt
\E_{-\infty}(Z)=0, \E_\infty(Z)=Z, \sum
_{n\in\mathbb{Z}} \|P_{n}Z \|_{2,\X}<\infty\biggr\} .
\end{equation}
It is not difficult to see that, setting $\|Z\|_{H_2}:=\sum_{n\in
\mathbb{Z}} \|P_{n}Z\|_{2,\X}$, $(H_2,\|\cdot\|_{H_2})$ is a Banach
space.

By our regularity conditions, we have, $Z=\sum_{k\in\mathbb{Z}}
P_kZ$ in $L^2(\Omega,\X)$ and $\P$-a.s.
Hence, writing $S_n =S_n(Z)=
\sum_{i=0}^{n-1} Z\circ\theta^i$, we have
\begin{eqnarray*}
S_n =\sum_{k\in\mathbb{Z}} \sum
_{i=0}^{n-1}(P_k Z) \circ
\theta^i .
\end{eqnarray*}
This splitting of $S_n$ into a series of martingales with (stationary)
increments has been used already in \cite{Wu} and \cite{CV}
in a similar context. This idea seems to appear first (explicitly)
in a paper by McLeish \cite{McLeish}.
We deduce that
\[
\M_2 (Z)\le\sum_{k\in\mathbb{Z}} \M_2
\bigl(P_k(Z)\bigr) .
\]
But, for every $k\in\mathbb{Z}$, $((P_kZ)\circ\theta^i)_{i\ge1}$
is a
stationary sequence of martingale differences. Hence, by Theorem~\ref{inemaxLIL}, for every $1\le p <2$, there exists $C_p$, such that
%
\begin{equation}
\label{maxhannan} \bigl\|\M_2 (Z)\bigr\|_{p,\infty}\le C_pD
\biggl(\sum_{k\in\mathbb{Z}} \|P_kZ\|
_{2,\X}\biggr) .
\end{equation}
%
We define a continuous operator $\D$ on $H_2$ with values in
$\{d \in L^2(\Omega,\F_1) \dvt  \E(d_1|\F_0)=0\}$, by setting, for
every $Z\in H_2$, $\D Z:= \sum_{n\in\mathbb{Z}} P_1 (Z\circ\theta^n)$.
Write $d=\D Z$. Let $M_n:=\sum_{i=0}^{n-1} d \circ\theta^i$. We want
to prove that
%
\begin{equation}
\label{es} |S_n-M_n|_\X=\mathrm{o}\bigl(\sqrt{nL
\bigl(L(n)\bigr)}\bigr)\qquad  \as
\end{equation}
Since $\M_2 (Z-d)\le\M_2(Z)+\M_2(d)$, using \eqref{maxhannan},
Theorem~\ref{inemaxLIL} and the Banach
principle (see the \hyperref[app]{Appendix}), we see that the set $\{Z\in H_2 \dvt  \mbox{\eqref{es} holds}\}$ is closed in $H_2$, and, by linearity, that set
is a
vector space.

Let $Z\in H_2$. Clearly, $Z=\sum_{k\in\Z}P_kZ$ in $H_2$. Hence it suffices
to prove \eqref{es} for $P_kZ$, for\vspace*{1pt} every $k\in\Z$. Now, $\D(P_k Z)
=(P_k Z)\circ\theta^{1-k}$. Let $k\le0$. We have
\begin{eqnarray*}
S_n(P_k Z)-M_n(P_k Z)&=& \sum
_{\ell=0}^{n-1} \bigl( (P_k Z)\circ
\theta ^\ell - (P_k Z)\circ\theta^{\ell+1-k} \bigr)
\\
&=& \sum_{\ell=0}^{-k} (P_k Z)\circ
\theta^\ell- \Biggl( \sum_{\ell=0}^{-k}
(P_k Z)\circ\theta^\ell \Biggr)\circ\theta^n
=\mathrm{o}(\sqrt n)\qquad  \as ,
\end{eqnarray*}
where we used that for any $X\in L^2(\Omega,\X)$,
$\sum_{n\ge1} \P(|X\circ\theta^n |_\X>\varepsilon\sqrt n)$, for
every $\varepsilon>0$, which implies that $X\circ\theta^n =\mathrm{o}(\sqrt n)$
\as, by the Borel--Cantelli lemma. The case $k\ge1$ may be handled
similarly.

\subsection{Proof of Theorem \texorpdfstring{\protect\ref{theoHannanp}}{2.8}}

As in the proof of Theorem~\ref{theoHannan}, we define a Banach space
\begin{eqnarray*}
\label{banhanp} H_p:=\biggl\{ Z\in L^p(\Omega, \X) \dvt
\E_{-\infty}(Z)=0, \E_\infty(Z)=Z, \|Z\|_{H_p}:=\sum
_{n\in\mathbb
{Z}} \|P_{n}Z \|_{p,\X}<\infty
\biggr\} .
\end{eqnarray*}

We see that
\begin{eqnarray*}
\label{maxhannanp} \|\M_p Z\|_{p,\infty}\le C_{p,r}D^{1/p}
\|Z\|_{H_p} ,
\end{eqnarray*}
where $C_{r,p}$ is the constant appearing in Proposition~\ref{mzmax},
and that the operator $\D$ may be extended
in a continuous operator from $H_p$ to $\{d\in L^p(\Omega,\F_1,\X) \dvt
\E_0(d)
=0\}$. Then, the proof\vadjust{\goodbreak} follows the one of
Theorem~\ref{theoHannan}. We first see that $|S_n-M_n|_\X=\mathrm{o}(n^{1/p})\ \as$ and then we use that the Marcinkiewicz--Zygmund strong law of
large number is known for $r$-smooth valued stationary martingale
differences, see, for example, \cite{Woyczinski}.

\subsection{Proof of Corollary \texorpdfstring{\protect\ref{corberg}}{2.12}}

We only have to prove that $\K_d$ is given as in the corollary.
By \eqref{Hannan}, we have $\sum_{n\in\mathbb{Z}}\|P_1X_n\|_{2,\X
}<\infty$.
Hence, for every $f,g\in\X^*$, we have, with absolute convergence of
all the
series,
\begin{eqnarray*}
\K_d(f,g)&=&\sum_{m,n\in\mathbb{Z}} \E
\bigl(P_1\bigl(f(X_n)\bigr)P_1
\bigl(g(X_m)\bigr)\bigr) =\sum_{m,n\in\mathbb{Z}} \E
\bigl(f(X_0)P_{1-n}\bigl(g(X_{m-n})\bigr)\bigr)
\\
&=&\sum_{m,n\in\mathbb{Z}} \E\bigl(f(X_0)P_{-n}
\bigl(g(X_{m})\bigr)\bigr) =\sum_{m\in\mathbb{Z}} \E
\bigl(f(X_0)g(X_{m})\bigr) .
\end{eqnarray*}

%

\section{Proof of Lemma \texorpdfstring{\protect\ref{Hanbis2}}{2.13}}

Since the sequences $(\|\E_{-n}(X)\|_{p,\H})$ and $(\|X-\E_{n}(X)\|
_{p,\H})$ are
non-increasing, \eqref{hanbis} is equivalent to
\[
\sum_{n\ge0}2^{n/2} \bigl\|\E_{-2^n}(X)
\bigr\|_{p,\H}<\infty \quad \mbox{and}\quad  \sum_{n\ge0}2^{n/2}
\bigl\|X-\E_{2^n}(X)\bigr\|_{p,\H}<\infty .
\]
In particular, $X$ is regular.

Assume $p=2$. For every $n\ge0$, using Cauchy--Schwarz and that $\E
(\langle P_{-k}X,P_{-\ell}X\rangle_\H)=0$ for
every $k\neq\ell$, we have
\begin{eqnarray*}
\Biggl(\sum_{k=2^n}^{2^{n+1}-1} \|P_{-k}X
\|_{2,\H} \Biggr)^2 \le2^n \sum
_{k\ge2^n} \E\bigl(|P_{-k}X|_\H^2
\bigr) \le2^n \E\bigl( \bigl|\E_{-2^n}(X)\bigr|_\H^2
\bigr) ,
\end{eqnarray*}
and
\begin{eqnarray*}
\Biggl(\sum_{k=2^n}^{2^{n+1}-1}
\|P_{k}X\|_{2,\H} \Biggr)^2 \le2^n
\sum_{k\ge2^n} \E\bigl(|P_{k}X|_\H^2
\bigr) \le2^n \E\bigl( \bigl|X-\E_{2^n}(X)\bigr|_\H^2
\bigr) .
\end{eqnarray*}

Assume $1<p<2$. By H\"older's inequality twice we have, with $1/p+1/q=1$,
\begin{eqnarray*}
\Biggl(\sum_{k=2^n}^{2^{n+1}-1} \|P_{-k}X
\|_{p,\H} \Biggr)^p &\le&2^{np/q} \E \Biggl(\sum
_{k= 2^n}^{2^{n+1}-1} |P_{-k}X|_\H^p
\Biggr)
\\
& \le&2^{np/2} \E \biggl(\biggl(\sum_{k\ge2^n}
|P_{-k}X|_\H^2\biggr)^{p/2} \biggr)
\le C2^{np/2} \bigl\|\E_{-2^n}(X)\bigr\|_{p,\H} ^p ,
\end{eqnarray*}
and
\begin{eqnarray*}
 \Biggl(\sum_{k=2^n}^{2^{n+1}-1}
\|P_{k}X\|_{p,\H} \Biggr)^2 &\le&2^{np/q}
\E \Biggl(\sum_{k= 2^n}^{2^{n+1}-1}
|P_{k}X|_\H^p \Biggr)
\\
& \le&2^{np/2} \E \biggl(\biggl(\sum_{k\ge2^n}
|P_{k}X|_\H^2\biggr)^{2/p} \biggr)
\le C2^{np/2} \bigl\|X-\E_{2^n}(X)\bigr\|_{p,\H}^p ,
\end{eqnarray*}
where we used Burkholder's inequality in Hilbert spaces, see \cite
{Burkholder}. Then, we conclude as above.

\subsection{Proof of Theorem \texorpdfstring{\protect\ref{dynsys}}{2.14}}

For every $n\ge0$ define $P^{(n)}:=\E^n-\E^{n+1}$. It suffices to prove
the theorem under the weaker condition
\[
\E^{\infty}(X)=0 \quad \mbox{and} \quad \sum_{n\ge0}
\bigl\|P^{(n)}(X)\bigr\|_{2,\H} <\infty .
\]
The fact that \eqref{conddynsys} implies the above condition may be
proved as
Lemma~\ref{Hanbis2}, using \eqref{relfro}.

Then, the proof may be done exactly as the proof of Theorem~\ref{theoHannan} except that we make use of reverse martingales.

\begin{appendix}\label{app}

\section{Proof of Proposition \texorpdfstring{\lowercase{\protect\ref{mzmax}}}{2.1}}\label{appA}
We start with the case
$d\in L^p(\Omega,\F_1,\P)$ and $\E_0(d)=0$. Define
\[
M^*=M^*(d):=\sup_{s\ge0}\frac{\max_{1\le n\le2^s}|S_n(d)|_\X}{2^{s/p}}.
\]
Let $s\ge0$. For every $2^s \le n\le2^{s+1}-1$, we have
\[
\frac{|S_n(d) |_\X}{n^{1/p}} \le\frac{\max_{1\le n\le2^s}
|S_n(d)|_\X}{2^{s/p}} \le M^* .
\]
Hence, it suffices to prove the result for $M^*$ instead of
$\M_p(d)$. Let $\lambda>0$.
We proceed by truncation. For every $s\ge0$, $k\ge1$ define
\begin{eqnarray*}
e_k^{(s)}&:=& d_k\mathbf{1}_{\{|d_k |_\X\le\lambda2^{s/p}
\}} ;\qquad
d_k^{(s)}:= e_k^{(s)} -\E
\bigl(e_{k}^{(s)}|\F_{k-1}\bigr) ;\\
  \tilde
e_k^{(s)}&:=& d_k-e_k^{(s)} ;\qquad
\tilde d_k^{(s)}:=d_k-d_k^{(s)}
;\\
  M_k^{(s)}&:=&\sum_{i=1}^k
d_i^{(s)} ; \qquad \tilde M_k^{(s)}:=
M_k-M_k^{(s)} .
\end{eqnarray*}
Let $\lambda>0$. Then,
\begin{eqnarray*}
&&\P\bigl( M^* >\lambda\bigr)
\\
&&\quad \le\sum_{s\ge0}\P\biggl(\frac{\max_{1\le n\le2^s}|\tilde
M_n^{(s)}|_\X}{2^{s/p}} >\lambda/2
\biggr) +\sum_{s\ge0}\P\biggl(\frac{\max_{1\le n\le
2^s}|M_n^{(s)}|_\X}{2^{s/p}} >
\lambda/2\biggr)
\\
&&\quad \le\frac{4}{\lambda}\sum_{s\ge0} 2^{(1-1/p)s}\E
\bigl(\bigl|\tilde e_1^{(s)}\bigr|_\X\bigr) +
\frac{2^{r}}{\lambda^r}\sum_{s\ge0}\frac{\E(
\max_{1\le n\le2^s}|M_n^{(s)}|_\X^r)}{2^{rs/p}} .
\end{eqnarray*}
Now, by Fubini and stationarity,
\[
\sum_{s\ge0} 2^{(1-1/p)s}\E\bigl(\bigl|\tilde
e_1^{(s)}\bigr|_\X\bigr) \le\frac{C\E
(|d_1|_\X^p)} {
\lambda^p}.
\]
To deal with the second term, we use Doob's maximal inequality in
$L^r$, for the submartingale $(|M_n|_\X)_{n\ge1}$, and \eqref{smooth}.
We obtain
\renewcommand{\theequation}{\arabic{equation}}
\setcounter{equation}{38}
\begin{eqnarray}
\label{doobr}\sum_{s\ge0}\frac{\E( \max_{1\le n\le
2^s}|M_n^{(s)}|_\X^r)}{2^{rs/p}} &\le&\sum
_{s\ge0}\frac{C_r}{2^{rs/p}\lambda^r}\E\bigl(\bigl|M_{2^s}^{(s)}\bigr|_\X
^r\bigr)\nonumber
\\[-8pt]\\[-8pt]
&\le& D^rC_r \sum_{s\ge0}2^{(1-r/p)s}
\E\bigl(\bigl|d_1^{(s)}\bigr|_\X^r\bigr) \le
\frac{D^rC_{r,p}\E(|d_1|_\X^p)} {
\lambda^{p-r}} ,\nonumber
\end{eqnarray}
which proves the proposition, in that case.
When $d\in L^2(\Omega,\F^{0},\P)$ and $\E^1(d)=0$, the proof is the same,
with the obvious changes, noticing that for every $n\ge1$,
$(S_n(d)-S_{n-k}(d))_{0\le k \le n}$ is a $(\F^{n-k})_{0\le k \le
n}$-martingale
and that $\max_{1\le k\le n}|S_k(d)|_\X\le2
\max_{1\le k\le n}|S_n(d)-S_{n-k}(d)|_\X$.

%




%


%



\section{Proof of Corollary \texorpdfstring{\lowercase{\protect\ref{flp}}}{3.4}}

Notice that, by \eqref{condcont0}, for every $x,h,h'\in\R$, we have
\setcounter{equation}{39}
\begin{equation}
\label{condcont} \bigl|f(x+h)-f\bigl(x+h'\bigr)\bigr|\le2^{r}\varphi
\bigl(\bigl|h-h'\bigr|\bigr) \bigl(1+|x|^r\bigr)+
2^{r-1}K\bigl(|h|^r+\bigl|h'\bigr|^r\bigr) .
\end{equation}
Recall that for every concave $\psi$ with $\psi(0)=0$,
$x\to\psi(x)/x$ is non-increasing on $]0,+\infty[$ and
$\psi$ is sub-additive.

We want to apply Theorem~\ref{theoHannan} and Lemma~\ref{Hanbis2}.
We shall evaluate $\|P_0(X_n)\|_2$, $\|\E_0(X_n)\|_2$
and $\|X_n-\E_n(X_n)\|$.

Enlarging our probability space if necessary, we assume that there exists
$(\xi'_n)$ an independent copy of $(\xi_n)$.

Then,
\begin{eqnarray*}
P_0X_n = \E_0 \bigl( f( A_n+h_n)
-f\bigl( A_n+h_n'\bigr) \bigr) ,
\end{eqnarray*}
where $A_n:=\sum_{k> -n} a_{-k}\xi_{n+k}'
+ \sum_{k>n} a_k \xi_{n-k}$, $h_n:= a_n\xi_0 $ and $h_n':= a_n\xi
_0' $.

In particular, we have,
by independence and using \eqref{condcont},
\[
\E\bigl((P_0X_n)^2\bigr) \le C_r
\bigl( \E\bigl(\varphi^2\bigl(|a_n|\bigl(|\xi_0|+\bigl|
\xi_0'\bigr|\bigr)\bigr)\bigr) \E\bigl(|A_n|^{2r}
\bigr) + |a_n|^{2r}\E\bigl(|\xi_0|^{2r}
\bigr) \bigr).
\]


We notice now that for every $\varphi\in\Lambda$, there exists $C>0$
such that, for every $n\ge1$
%
\begin{equation}
\label{techconc} \E\bigl(\varphi^2\bigl(|a_n|\bigl(|
\xi_0|+\bigl|\xi_0'\bigr|\bigr)\bigr)\bigr)\le C
\varphi^2\bigl(|a_n|\bigr) .
\end{equation}
This follows from Jensen's inequality and the sub-additivity of
$\varphi^2$
(using that $\xi_0\in L^1(\Omega,\F,\P)$) when $\varphi^2$ is
sub-additive, and it is obvious when $\varphi(x)=\min(1,x^\alpha)$
(using that $\xi_0\in L^{2\alpha}(\Omega,\allowbreak \F,\P)$).

Clearly, $\E(|A_n|^{2r}) \le (\sum_{k\in\mathbb{Z}}|a_k|\|\xi
_0\|_{2r} )^{2r}$.

Since $x\to\varphi^2(x)/x$ is non-increasing, when $\varphi^2$ is
concave, we see that whenever
$\varphi\in\Lambda$, $|a_n|^{2r}\le C \varphi^2(|a_n|)$.

This finishes the proof of Corollary~\ref{flp} under the assumption
on $P_0(X_n)$.

We shall now evaluate $\|\E_0(X_n)\|_2$, the case of $\|X_n-\E
_n(X_n)\|_2$
may be treated similarly. We have
\[
\E_0(X_n)= \E_0\bigl(f(B_n+k_n)-f
\bigl(B_n-k_n'\bigr)\bigr) ,
\]
where $B_n:=\sum_{k> -n} a_{-k}\xi_{n+k}$, $k_n=\sum_{k\ge n} a_k
\xi_{n-k}$ and $k_n'=\sum_{k\ge n} a_k \xi_{n-k}'$.
Hence, using \eqref{condcont},
\[
\|E_0(X_n)\|_2^2 \le
C_r \bigl( \E\bigl(\varphi^2\bigl(|k_n|+|k_n'|
\bigr) \E\bigl(|A_n|^{2r}\bigr) + 2\|k_n
\|_{2r}^{2r} \bigr) \bigr) .
\]
When $\varphi^2$ is concave, by Jensen's inequality,
\[
\E\bigl(\varphi^2\bigl(|k_n|+\bigl|k_n'\bigr|
\bigr)\bigr)\le\varphi^2\biggl(2\E\bigl(|\xi_0|\bigr)\sum
_{k\ge
n}|a_k|\biggr)\le\bigl(1+2\E\bigl(|
\xi_0|\bigr)\bigr)\varphi^2\biggl(\sum
_{k\ge n}|a_k|\biggr).
\]
When $\varphi(x)=\min(1,x^\alpha)$, assuming that $1/2\le\alpha\le
1$ (otherwise we are in the previous case),
we have
\[
\E\bigl(\varphi^2\bigl(|k_n|+\bigl|k_n'\bigr|
\bigr)\bigr)\le \biggl(\sum_{k\ge n}|a_k|\|
\xi_0\| _{2\alpha} \biggr)^{2\alpha}\le C
\varphi^2\biggl(\sum_{k\ge n}|a_k|
\biggr).
\]
Clearly, $\E(|B_n|^{2r}) \le (\sum_{k\in\mathbb{Z}}|a_k|\|\xi
_0\|_{2r} )^{2r}$.

Finally, we have
\[
\|k_n\|_{2r}^{2r}\le\|\xi_0
\|_{2r}^{2r} \biggl( \sum_{k\ge n}
|a_k| \biggr)^{2r} .
\]
Since $x\to\varphi^2(x)/x$ is non-decreasing, when $\varphi^2$ is
concave, we see that whenever
$\varphi\in\Lambda$,
\[
\|k_n\|_{2r}^{2r}\le C \varphi^2
\biggl(\sum_{k\ge n} |a_k|\biggr) .
\]

\section{The Banach principle}

The following is an extension of the Banach principle as
stated in Theorem~7.2, page 64 of \cite{Krengel}.

\begin{prop}\label{banprin}
Let $(\Omega,\F,\P)$ be a probability space and $\X, \B$ be
Banach spaces. Let $\C$ be a vector space of measurable
functions from $\Omega$ to $\X$. Let $(T_n)_{n\ge1}$ be a sequence
of linear maps from $\B$ to $\C$. Assume that there exists
a positive decreasing function $L$ on $]0,+\infty[$, with
$\lim_{\lambda\to\infty}L(\lambda)=0$, such that
\setcounter{equation}{41}
\begin{equation}
\label{bancond} \P\Bigl(\sup_{n\ge1} |T_n
x|_\X> \lambda|x|_\B\Bigr)\le L(\lambda)\qquad  \forall
\lambda>0,x\in\B .
\end{equation}
Then the set $\{x\in\B \dvt  (T_n x)_{n\ge1} \mbox{ is \as\ relatively compact
in $\X$}\}$
and the set $\{x\in\B\dvt  |T_n x|_\X\to0 \ \as\}$
are closed in $\B$.
\end{prop}

\begin{pf} We prove that the first set is closed, the proof for
the second one being similar, but easier.
Let $x\in\B$ and $(x_m)_{m\ge1}\subset\B$ be such that
$|x_m-x|_\B\displaystyle \mathop{\longrightarrow}_{m\to\infty}0$ and such that for
every $m\ge1$, $(T_n x_m)_{n\ge1}$ is \as\ relatively compact in
$\X$. We want to prove that $(T_nx)_{n\ge1}$ is $\P$-a.s. relatively compact.

By \eqref{bancond}, for every integers $m,p\ge1$ (assume that $x\neq x_m$
otherwise there is nothing to do)
\begin{eqnarray*}
\P\Bigl(\sup_{n\ge1} \bigl|T_n (x-x_m)\bigr|_\X>
1/p\Bigr)\le L \biggl(\frac{1}{p|x-x_m|_\B} \biggr)\qquad  \forall\lambda>0,x\in\B .
\end{eqnarray*}
Since $\lim_{\lambda\to\infty}L(\lambda)=0$, there exists a
subsequence $(m_k)_{k\ge1}$ and a set
$\Omega_0\in\F$ with $\P(\Omega_0)=1$, such that for every $\omega
\in\Omega_0$,
\[
\sup_{n\ge1} \bigl|T_n (x-x_{m_k})\bigr|_\X(
\omega)\mathop{\longrightarrow}_{k\to\infty
}0 .
\]

There exists $\Omega_1\in\F$, with $\P(\Omega_1)=1$,
such that, for every $\omega\in\Omega_1$ and every $k\ge1$, $((T_n
x_{m_k})(\omega))_{n\ge1}$
is relatively compact in $\X$.

Let $\omega\in\Omega_0\cap\Omega_1$ be fixed. Let $\varphi_0$ be
an increasing function from $\mathbb{N}$ to
$\mathbb{N}$. We want to prove that $(T_{\varphi_0(n)}x(\omega
))_{n\ge1}$ admits a convergent subsequence.

For every $k\ge1$, $((T_{\varphi_0(n)} x_{m_k})(\omega))_{n\ge1}$
admits a Cauchy subsequence.
We construct by induction some
increasing functions $(\varphi_k)_{k\ge1}$ such that, for every
$k\ge1$, setting $\psi_k:=\varphi_0\circ\varphi_1 \circ\cdots
\circ\varphi_k$, we have for every $p\ge n\ge1$,
\[
\bigl|T_{\psi_k(n)}x_{m_k}(\omega) -T_{\psi_k(p)}x_{m_k}(
\omega)\bigr|_\X\le 1/n .
\]
Then, $(T_{\psi_n(n)}x(\omega))$ is Cauchy. Indeed, for every $N\ge
1$, and every $p> n \ge N$, we have
\begin{eqnarray*}
&&\bigl|T_{\psi_n(n)}x(\omega) -T_{\psi_p(p)}x(\omega)\bigr|_\X
\\
&&\quad \le\bigl|T_{\psi
_n(n)}x_{m_n}(\omega) -T_{(\psi_n\circ
\varphi_{n+1}\circ\cdots\circ\varphi_p)(p)}x_{m_n}(
\omega)\bigr|_\X +2 \sup_{r\ge1} \bigl|T_r(x_{m_n}-x)\bigr|_\X
\mathop{\longrightarrow}_{N\to\infty}0 ,
\end{eqnarray*}
and the result follows. \end{pf}

%


\section{Identification of the cluster set}\label{cluster}

Denote by $\X^*$ the topological dual of $\X$. Let $X \in L^2(\Omega
,\X)$ such that $\E(X)=0$. Following Kuelbs \cite{Kuelbs}
(we refer to \cite{Kuelbs} for more details on the construction below),
we define a bounded linear
operator $\S$ from $\X^*$ to $\X$ and a bounded \emph{symmetric}
bilinear operator $\K$ from $\X^*\times\X^*$
to $\mathbb{R}$, by
\begin{eqnarray*}
\S\bigl(x^*\bigr)&=& \E\bigl(x^*(X)X\bigr)\qquad  \forall x^*\in\X^* ,
\\
\K\bigl(x^*,y^*\bigr)& =&\E\bigl(x^*(X)y^*(X)\bigr) =y^*\bigl(\S\bigl(x^*\bigr)
\bigr)=x^*\bigl(\S\bigl(y^*\bigr)\bigr)\qquad  \forall x^*,y^*\in\X^* .
\end{eqnarray*}

Let $\H_X$ be the closure of the range of $\S$ with respect to the
following inner
product:
\[
\bigl\langle\S x^*,\S y^* \bigr\rangle_{\H_X} = \K\bigl(x^*,y^*\bigr) .
\]
Notice that the definition of $\langle\cdot, \cdot\rangle_{\H_X}$
does not depend on the chosen representatives (i.e., if $x^*\in\Ker \S$,
$\langle\S x^*,\S y^* \rangle_{\H_X}=0$ for every $y^*\in\X^*$)
and that this inner product is really \emph{positive
definite}.

Finally, denote by $K=K_X$, the unit ball of $(\H_X,\|\cdot\|_{\H
_X})$, $K$ is compact by (iv), Lemma~2.1 of \cite{Kuelbs}. We recall
an important result of Kuelbs, see
\cite{Kuelbs}, Theorem~3.1, II, where we denote by $C(\{x_n\})$ the
cluster set of a sequence $(x_n)\subset\X$.

\begin{prop}[(Kuelbs \cite{Kuelbs})]\label{Kuelbs}
Let $X\in L^2(\Omega,\X)$. Assume that $(X_n)_{n\ge0}$ satisfies the
{CLIL} and that,
\setcounter{equation}{42}
\begin{equation}
\label{weakLIL1} \limsup_n \frac{S_n(x^*(X))}{\sqrt{2nL(L(n))}} = \bigl\| x^*(X)
\bigr\|_2\qquad  \as \qquad \forall x^*\in\X^* .
\end{equation}
Then,
%
\begin{eqnarray}
\label{comp0}C \biggl(\biggl\{\frac{S_n(X)}{\sqrt{2nL(L(n))}} \biggr\} \biggr)=K \qquad \as,
\end{eqnarray}
 and
 \begin{eqnarray}
\label{norm0} \limsup_{n} \frac{|S_n(X)|_\X}{\sqrt {2nL(L(n))} }=
\sup_{x^*\in\X^*,|x^*|_{\X^*}\le1} \bigl\|x^*(X)\bigr\|_2 \le \|X\|_{2,\X}
\qquad \as
\end{eqnarray}
\end{prop}

\end{appendix}

\section*{Acknowledgements}
The author is very grateful to the referees. One of them provided
a very detailed and accurate report that surely helped to improve the
presentation
of the paper. The author is also grateful to Florence Merlev\`ede for
reading a first version of the paper and for her comments.



\printhistory

\end{document}